\documentclass[11pt]{article}

\usepackage[a4paper]{geometry}
\usepackage{indentfirst}
\usepackage{hyperref}
\usepackage{breakurl}
\usepackage{amssymb, amsmath}

\hypersetup{pdfstartview=FitH, pdfpagemode=UseNone}

\newcounter{lemma}
\newenvironment{lemma}{\refstepcounter{lemma} \emph{Lemma \thelemma.}}{\ignorespacesafterend}
\newcounter{theorem}
\newenvironment{theorem}{\refstepcounter{theorem} \emph{Theorem \thetheorem.}}{\ignorespacesafterend}
\newenvironment{proof}{\emph{Proof.}}{$\square$}
\newenvironment{proof*}[1]{\emph{#1.}}{$\square$}

\allowdisplaybreaks

\begin{document}

\title{\textbf{One Curious Identity\\ Counting Graceful Labelings}}
\author{Nikolai Beluhov}
\date{}

\maketitle

\begin{center} \parbox{352pt}{\setlength{\parindent}{12pt} \footnotesize \emph{Abstract}. Let $a$ and $b$ be positive integers with prime factorisations $a = p_1^np_2^n$ and $b = q_1^nq_2^n$. We prove that the number of essentially distinct $\alpha$-graceful labelings of the complete bipartite graph $K_{a, b}$ equals the alternating sum of fourth powers of binomial coefficients $(-1)^n[\binom{2n}{0}^4 - \binom{2n}{1}^4 + \binom{2n}{2}^4 - \binom{2n}{3}^4 + \cdots + \binom{2n}{2n}^4]$.} \end{center}

\section{Introduction} \label{intro}

Consider a bipartite graph $G$ with $a$ vertices in its left part, $b$ vertices in its right part, and $e$ edges.

Let us label each vertex of $G$ with one of the nonnegative integers $0$, $1$, \ldots, $e$ so that every label in the left part of $G$ is strictly smaller than every label in the right part of $G$. Then let us also label every edge of $G$ with the difference between the labels of its endpoints. When the labels on the edges of $G$ turn out to be the positive integers $1$, $2$, \ldots, $e$, each one occurring exactly once, we say that our labeling is \emph{$\alpha$-graceful}. These labelings were introduced by Alexander Rosa in \cite{R}.

We consider two $\alpha$-graceful labelings to be essentially the same when they can be transformed into one another by automorphisms of $G$.

Let $\mathcal{A}(a, b)$ be the number of essentially distinct $\alpha$-graceful labelings of the complete bipartite graph $K_{a, b}$. Equivalently, $\mathcal{A}(a, b)$ is also the number of ordered pairs of sets $S'$ and $S''$ of nonnegative integers such that $|S'| = a$, $|S''| = b$, both of $S'$ and $S''$ are subsets of $\{0, 1, \ldots, ab\}$, and $S'' - S' = \{y - x \mid x \in S' \text{ and } y \in S''\} = \{1, 2, \ldots, ab\}$.

Donald Knuth studies $\mathcal{A}(a, b)$ in \cite{K}. We give a brief overview of some results obtained there in Section \ref{prelim}.

Our goal in this short note will be to prove the following theorem.

\medskip

\begin{theorem} \label{main} Suppose that the prime factorisations of $a$ and $b$ are $a = p_1^np_2^n$ and $b = q_1^nq_2^n$. Then \begin{align*} \mathcal{A}(a, b) &= \sum_{i = 0}^{2n} (-1)^{n + i}\binom{2n}{i}^4\\ &= (-1)^n\left[\binom{2n}{0}^4 - \binom{2n}{1}^4 + \binom{2n}{2}^4 - \binom{2n}{3}^4 + \cdots + \binom{2n}{2n}^4\right]. \end{align*} \end{theorem}

\medskip

Note that we do not require the sets of prime factors $\{p_1, p_2\}$ of $a$ and $\{q_1, q_2\}$ of $b$ to be disjoint.

Nicolaas de Bruijn studied the asymptotic behaviour of the alternating sums of like powers of binomial coefficients \begin{align*} S(k, n) &= \sum_{i = 0}^{2n} (-1)^{n + i}\binom{2n}{i}^k\\ &= (-1)^n\left[\binom{2n}{0}^k - \binom{2n}{1}^k + \binom{2n}{2}^k - \binom{2n}{3}^k + \cdots + \binom{2n}{2n}^k\right] \end{align*} in \cite{B}. Thus these sums are now known as \emph{de Bruijn's $S(k, n)$}. With this notation, we can restate Theorem \ref{main} quite succinctly as $\mathcal{A}(p_1^np_2^n, q_1^nq_2^n) = S(4, n)$.

The author discovered the identity of Theorem \ref{main} with the help of the On-Line Encyclopedia of Integer Sequences. The results of Section \ref{prelim} allow us to determine $\mathcal{A}(p_1^np_2^n, q_1^nq_2^n)$ experimentally for small $n$. When we enter the first three values $14$, $786$, and $61340$ into the OEIS search engine, immediately we are presented with the OEIS entry \cite{OEIS} for de Bruijn's $S(4, n)$.

\section{Preliminaries} \label{prelim}

The results summarised in this section are all established in \cite{K}. The wording there is somewhat different; the one we prefer here, in terms of transformations, is also due to Knuth.

We can construct larger $\alpha$-graceful labelings out of smaller ones as follows. Take any $\alpha$-graceful labeling of $K_{a, b}$ with labels $S'$ on the left and $S''$ on the right, and let $c$ be any positive integer with $c \ge 2$. Then the sets of labels $T' = cS' + \{0, 1, \ldots, c - 1\} = \{cx + z \mid x \in S' \text{ and } 0 \le z < c\}$ and $T'' = cS'' = \{cy \mid y \in S''\}$ determine an $\alpha$-graceful labeling of $K_{ac, b}$. We call this transformation \emph{multiplication by $c$ on the left}.

Similarly, \emph{multiplication by $c$ on the right} produces an $\alpha$-graceful labeling of $K_{a, bc}$ with label sets $U' = cS'$ on the left and $U'' = cS'' - \{0, 1, \ldots, c - 1\}$ on the right.

Every $\alpha$-graceful labeling of $K_{a, b}$ can be obtained from the unique $\alpha$-graceful labeling of $K_{1, 1}$ by means of some series $\mathcal{S}$ of multiplications on the left and right.

We say that $\mathcal{S}$ is \emph{canonical} when it alternates between multiplications on the left and multiplications on the right.

Observe that two successive multiplications on the left by $c'$ and then $c''$ yield the same net result as one single multiplication on the left by $c'c''$. This allows us to collapse any subseries of $\mathcal{S}$ which consists entirely of multiplications on the left into one single multiplication on the left. The same things hold true of multiplication on the right as well. Therefore, we can assume without loss of generality that $\mathcal{S}$ is canonical.

Conversely, distinct canonical series of multiplications yield distinct $\alpha$-graceful labelings. Therefore, $\mathcal{A}(a, b)$ equals the number of distinct canonical series of multiplications on the left and right which transform $K_{1, 1}$ into $K_{a, b}$.

Let $a = p_1^{a_1}p_2^{a_2} \cdots p_r^{a_r}$ and $b = q_1^{b_1}q_2^{b_2} \cdots q_s^{b_s}$ be the prime factorisations of $a$ and $b$. By the preceding discussion, $\mathcal{A}(a, b)$ depends only on the unordered pair of multisets $\{a_1, a_2, \ldots, a_r\}$ and $\{b_1, b_2, \ldots, b_s\}$. Thus we may just as well do away with the prime factors altogether and write simply $\mathcal{A}(\{a_1, a_2, \ldots, a_r\}, \{b_1, b_2, \ldots, b_s\})$. With this notation, the identity of Theorem \ref{main} becomes $\mathcal{A}(\{n, n\}, \{n, n\}) = S(4, n)$.

\section{One General Theorem} \label{thm}

Let $P$ be a polynomial of degree $d$ with real coefficients. Then there are unique real numbers $\mathcal{C}(0, P)$, $\mathcal{C}(1, P)$, \ldots, $\mathcal{C}(d, P)$ such that \[P(x) = \sum_{i = 0}^d \mathcal{C}(i, P)\binom{x}{i}.\]

We define also $\mathcal{C}(k, P) = 0$ for all $k > d$.

The theory of finite difference allows us to express the $\mathcal{C}(k, P)$ neatly in terms of the values of $P$. They are given by \[\mathcal{C}(k, P) = \sum_i (-1)^{k + i}\binom{k}{i}P(i). \tag{$\ast$}\]

Here and henceforth, we write ``$\sum_i$'' as shorthand for ``$\sum_{i = 0}^\infty$'', provided that the summand vanishes for all sufficiently large $i$ and so the sum actually contains only finitely many nonzero terms.

Let $Q$ be a polynomial with real coefficients as well. We define \[P \bigcirc Q = \sum_i \mathcal{C}(i, P)\mathcal{C}(i, Q).\]

Thus $\bigcirc$ assigns a real number to every pair of polynomials with real coefficients.

Observe that \[\begin{gathered} P \bigcirc Q = Q \bigcirc P,\\ (c'P' + c''P'') \bigcirc Q = c'(P' \bigcirc Q) + c''(P'' \bigcirc Q), \end{gathered}\] and \[P \bigcirc (d'Q' + d''Q'') = d'(P \bigcirc Q') + d''(P \bigcirc Q'').\]

That is, $\bigcirc$ is commutative and linear in both arguments.

Furthermore, \[\binom{x}{k} \bigcirc P(x) = \mathcal{C}(k, P)\] for all nonnegative integers $k$.

Given a positive integer $a$ with prime factorisation $a = p_1^{a_1}p_2^{a_2} \cdots p_r^{a_r}$, we define \[\mathcal{F}_a(x) = \prod_{i = 1}^r \binom{x + a_i}{a_i}.\]

Thus $\mathcal{F}_a(x)$ is a polynomial with real coefficients of degree $a_1 + a_2 + \cdots + a_r$.

To calculate $\mathcal{A}(a, b)$, we rely on the following theorem.

\medskip

\begin{theorem} \label{aux} The number of essentially distinct $\alpha$-graceful labelings of the complete bipartite graph $K_{a, b}$ is given by \[\mathcal{A}(a, b) = [\mathcal{F}_a(x - 1) \bigcirc \mathcal{F}_b(x)] + [\mathcal{F}_a(x) \bigcirc \mathcal{F}_b(x - 1)]\] in all cases except for $a = b = 1$, when $\mathcal{A}(1, 1) = 1$. \end{theorem} 

\medskip

In the exceptional case of $a = b = 1$, the formula of Theorem \ref{aux} yields $2$ instead of the correct answer $1$. The reasons for this discrepancy will become clear a little bit later on.

For the proof, first of all we define $d_0$, $d_1$, \ldots, $d_k$ to be a \emph{partial sequence for $a$ with $k$ steps} when $1 = d_0 < d_1 < \cdots < d_k = a$ and $d_i$ divides $d_{i + 1}$ for all $i$. Let $\vartheta(k, a)$ be the number of distinct partial sequences for $a$ with $k$ steps.

\medskip

\begin{lemma} \label{part} The number of distinct partial sequences for $a$ with $k$ steps is given by \[\vartheta(k, a) = \mathcal{C}(k, \mathcal{F}_a(x - 1)).\] \end{lemma} 

\medskip

Of course, then also the number of distinct partial sequences for $b$ with $\ell$ steps will be given by $\vartheta(\ell, b) = \mathcal{C}(\ell, \mathcal{F}_b(x - 1))$.

\medskip

\begin{proof} Observe that $\vartheta$ satisfies the recurrence relation \[\vartheta(k + 1, a) = \sum_{d \mid a \text{ and } d \neq a} \vartheta(k, d)\] as well as the initial conditions $\vartheta(0, 1) = 1$ and $\vartheta(0, a) = 0$ for all $a \ge 2$. Let us now see that $\mathcal{C}(k, \mathcal{F}_a(x - 1))$ satisfies them as well.

The initial conditions are straightforward enough. When $a = 1$, $\mathcal{F}_a(x)$ becomes the empty product and so $\mathcal{C}(0, \mathcal{F}_a(x - 1)) = \mathcal{C}(0, 1) = 1$, too. Otherwise, when $a \ge 2$, we get $\mathcal{F}_a(-1) = 0$ because all terms in the product vanish. Then $\mathcal{C}(0, \mathcal{F}_a(x - 1)) = \mathcal{F}_a(0 - 1) = 0$ as well.

We go on to the recurrence relation. We must check that \[\mathcal{C}(k + 1, \mathcal{F}_a(x - 1)) = \sum_{d \mid a \text{ and } d \neq a} \mathcal{C}(k, \mathcal{F}_d(x - 1)).\]

To begin with, let us evaluate the simpler sum $\sum_{d \mid a} \mathcal{C}(k, \mathcal{F}_d(x - 1))$. Once we are done, we will subtract out from the total the superfluous term $\mathcal{C}(k, \mathcal{F}_a(x - 1))$.

Since \[\sum_{d \mid a} \mathcal{C}(k, \mathcal{F}_d(x - 1)) = \mathcal{C}\left(k, \sum_{d \mid a} \mathcal{F}_d(x - 1)\right),\] we focus on the polynomial $\sum_{d \mid a} \mathcal{F}_d(x - 1)$.

Given any divisor $d$ of $a$, let $d = p_1^{\delta_1}p_2^{\delta_2} \cdots p_r^{\delta_r}$ so that all of the $\delta_i$ are nonnegative integers. Note that this is not quite the prime factorisation of $d$ since we allow some or all of the $\delta_i$ to be zeroes. Note also that still \[\mathcal{F}_d(x) = \prod_{i = 1}^r \binom{x + \delta_i}{\delta_i},\] because each vanishing $\delta_i$ contributes only a unit multiplicand to the product.

Then \begin{align*} \sum_{d \mid a} \mathcal{F}_d(x - 1) &= \sum_{d \mid a} \prod_{i = 1}^r \binom{x - 1 + \delta_i}{\delta_i}\\ &= \prod_{i = 1}^r \sum_{\delta_i = 0}^{a_i} \binom{x - 1 + \delta_i}{\delta_i}. \end{align*}

To simplify the latter expression, we require the hockey-stick identity \[\sum_{i = 0}^m \binom{n + i}{n} = \binom{m + n + 1}{n + 1}.\]

For this and all other well-known binomial coefficient identities that we cite, we use Henry Gould's \cite{G} as a standard reference on the subject. We find the hockey-stick identity in it under number (1.52).

Thus \[\sum_{\delta_i = 0}^{a_i} \binom{n - 1 + \delta_i}{\delta_i} = \sum_{\delta_i = 0}^{a_i} \binom{n - 1 + \delta_i}{n - 1} = \binom{n + a_i}{n} = \binom{n + a_i}{a_i}\] for all positive integers $n$.

Consequently, quite simply \begin{align*} \sum_{d \mid a} \mathcal{F}_d(x - 1) &= \prod_{i = 1}^r \binom{x + a_i}{a_i}\\ &= \mathcal{F}_a(x). \end{align*}

Therefore, \[\sum_{d \mid a} \mathcal{C}\left(k, \mathcal{F}_d(x - 1)\right) = \mathcal{C}(k, \mathcal{F}_a).\]

Let $P$ be any polynomial with real coefficients. From $P(x) = \sum_i \mathcal{C}(i, P)\binom{x}{i}$, we derive \begin{align*} P(x + 1) &= \sum_i \mathcal{C}(i, P)\binom{x + 1}{i}\\ &= \mathcal{C}(0, P) + \sum_{i = 1}^\infty \mathcal{C}(i, P)\left[\binom{x}{i - 1} + \binom{x}{i}\right]\\ &= \sum_i \left[\mathcal{C}(i, P) + \mathcal{C}(i + 1, P)\right]\binom{x}{i}. \end{align*}

Consequently, \[\mathcal{C}(k, P(x + 1)) = \mathcal{C}(k, P) + \mathcal{C}(k + 1, P).\]

With $P(x) = \mathcal{F}_a(x - 1)$, we finally conclude that \begin{align*} \sum_{d \mid a \text{ and } d \neq a} \mathcal{C}(k, \mathcal{F}_d(x - 1)) &= \left[\sum_{d \mid a} \mathcal{C}(k, \mathcal{F}_d(x - 1))\right] - \mathcal{C}(k, \mathcal{F}_a(x - 1))\\ &= \mathcal{C}(k, \mathcal{F}_a) - \mathcal{C}(k, \mathcal{F}_a(x - 1))\\ &= \mathcal{C}(k + 1, \mathcal{F}_a(x - 1)), \end{align*} as needed. \end{proof}

\medskip

We are ready to tackle Theorem \ref{aux}.

\medskip

\begin{proof*}{Proof of Theorem \ref{aux}} Let $\mathcal{S}$ be one canonical series of multiplications transforming $K_{1, 1}$ into $K_{a, b}$. Suppose, to begin with, that $\mathcal{S}$ starts with a multiplication on the left and ends with a multiplication on the right. In a moment we will account for the other three cases as well.

Consider the complete bipartite graphs into which $\mathcal{S}$ successively transforms $K_{1, 1}$. Let the entire sequence of such graphs be \[K_{1, 1} = K_{d_0, e_0} \to K_{d_1, e_0} \to K_{d_1, e_1} \to K_{d_2, e_1} \to K_{d_2, e_2} \to \cdots \to K_{d_k, e_k} = K_{a, b}.\] Thus the multiplications of $\mathcal{S}$ are, in this order, by $d_1$ on the left, by $e_1$ on the right, by $d_2/d_1$ on the left, by $e_2/e_1$ on the right, and so on. Then $d_0$, $d_1$, \ldots, $d_k$ is a partial sequence for $a$ and $e_0$, $e_1$, \ldots, $e_k$ is a partial sequence for $b$.

We can similarly decompose every canonical series $\mathcal{S}$ into a pair of partial sequences for $a$ and $b$. When $\mathcal{S}$ starts and ends with multiplications on opposite sides, the two partial sequences will have the same number of steps, as above. Otherwise, when $\mathcal{S}$ starts and ends with multiplications on the same side, the partial sequence which corresponds to that side will be one step longer.

Conversely, consider a pair of partial sequences $D$ for $a$ and $E$ for $b$ with $k$ and $\ell$ steps, respectively, such that $k$ and $\ell$ differ by at most one. Then we can interleave $D$ and $E$ so as to obtain a canonical series $\mathcal{S}$ with $k$ multiplications on the left and $\ell$ multiplications on the right.

When $|k - \ell| = 1$, the interleaving can only happen in one unique way, with the longer partial sequence corresponding to the first and last multiplications of $\mathcal{S}$. Otherwise, when $k = \ell \ge 1$, the interleaving can happen in two distinct ways. Finally, when $k = \ell = 0$, both of $D$ and $E$ become the empty sequence, and so once again they can only be interleaved uniquely to yield the empty series of multiplications which transforms $K_{1, 1}$ into $K_{1, 1}$.

Therefore, \begin{align*} \mathcal{A}(a, b) &= \sum_k \vartheta(k, a)\vartheta(k + 1, b) + \sum_\ell \vartheta(\ell + 1, a)\vartheta(\ell, b) + {}\\ &\phantom{{} = {}} \vartheta(0, a)\vartheta(0, b) + 2\sum_{i = 1}^\infty \vartheta(i, a)\vartheta(i, b). \end{align*}

When at least one of $a$ and $b$ is greater than or equal to two, we get $\vartheta(0, a)\vartheta(0, b) = 0$ and so we can rewrite the above as \[\mathcal{A}(a, b) = \sum_i \vartheta(i, a)[\vartheta(i, b) + \vartheta(i + 1, b)] + \sum_i \vartheta(i, b)[\vartheta(i, a) + \vartheta(i + 1, a)].\]

The same trick does not work when $a = b = 1$, because then $\vartheta(0, a)\vartheta(0, b) = 1$ does not vanish. This is why Theorem \ref{aux} fails in that special case. Suppose, for now on, that indeed at least one of $a$ and $b$ is greater than or equal to two.

By Lemma \ref{part}, $\vartheta(i, a) = \mathcal{C}(i, \mathcal{F}_a(x - 1))$, and similarly for $\vartheta(i + 1, a)$, $\vartheta(i, b)$, and $\vartheta(i + 1, b)$. Then, as in the proof of Lemma \ref{part}, \begin{align*} \vartheta(k, a) + \vartheta(k + 1, a) &= \mathcal{C}(k, \mathcal{F}_a(x - 1)) + \mathcal{C}(k + 1, \mathcal{F}_a(x - 1))\\ &= \mathcal{C}(k, \mathcal{F}_a), \end{align*} and similarly for $\vartheta(\ell, b) + \vartheta(\ell + 1, b)$.

Therefore, in the end we arrive at \begin{align*} \mathcal{A}(a, b) &= \sum_i \mathcal{C}(i, \mathcal{F}_a(x - 1))\mathcal{C}(i, \mathcal{F}_b) + \sum_i \mathcal{C}(i, \mathcal{F}_a)\mathcal{C}(i, \mathcal{F}_b(x - 1))\\ &= [\mathcal{F}_a(x - 1) \bigcirc \mathcal{F}_b(x)] + [\mathcal{F}_a(x) \bigcirc \mathcal{F}_b(x - 1)], \end{align*} as needed. \end{proof*}

\section{One Example} \label{exp} 

The material in this section is not part of our proof of Theorem \ref{main}. (With the exception of the first paragraph of the proof of Lemma \ref{exl}. The point it makes is crucial also for our proof of Lemma \ref{square}.) However, we include it anyway because it is a good example of how Theorem \ref{aux} works in practice.

The treatment of $\mathcal{A}(a, b)$ in \cite{K} contains also the following result.

\medskip

\begin{theorem} \label{ext} Suppose that the prime factorisations of $a$ and $b$ are $a = p^n$ and $b = q_1^{b_1}q_2^{b_2} \cdots q_s^{b_s}$. Then \[\mathcal{A}(a, b) = \prod_{i = 1}^s \binom{n + b_i}{b_i}.\] \end{theorem} 

\medskip

The proof in \cite{K} uses the theory of traces. Here, we give a different proof which employs Theorem \ref{aux} instead.

By Theorem \ref{aux}, \[\mathcal{A}(a, b) = \left[\binom{x + n - 1}{n} \bigcirc \mathcal{F}_b(x)\right] + \left[\binom{x + n}{n} \bigcirc \mathcal{F}_b(x - 1)\right].\]

To evaluate the right-hand side, first we learn to express $\binom{x + m}{n} \bigcirc P(x)$ in terms of the values of $P$.

\medskip

\begin{lemma} \label{exl} For all nonnegative integers $m$ with $m \le n$ and all polynomials $P$ with real coefficients, \[\binom{x + m}{n} \bigcirc P(x) = \sum_i (-1)^i\binom{n - m}{i}P(n - i).\] \end{lemma} 

\medskip

\begin{proof} Observe that it suffices to consider only the case when $P(x)$ is of the form $\binom{x}{k}$. Since $\bigcirc$ is linear in both arguments, then it would follow immediately that Lemma \ref{exl} holds also for all linear combinations of polynomials of that form. But we already know that in fact every polynomial with real coefficients is such a linear combination.

With $P(x) = \binom{x}{k}$, we get \[\binom{x}{k} \bigcirc \binom{x + m}{n} = \mathcal{C}\left(k, \binom{x + m}{n}\right).\]

On the other hand, by the well-known binomial coefficient identity \[\binom{\ell + m}{n} = \sum_i \binom{\ell}{i}\binom{m}{n - i}\] (number (3.1) in \cite{G}) we conclude that also \[\binom{x + m}{n} = \sum_i \binom{m}{n - i}\binom{x}{i}.\]

Therefore, quite simply \[\mathcal{C}\left(k, \binom{x + m}{n}\right) = \binom{m}{n - k}.\]

With this, Lemma \ref{exl} boils down to \[\binom{m}{n - k} = \sum_i (-1)^i\binom{n - m}{i}\binom{n - i}{k},\] and this is yet another well-known binomial coefficient identity. We find it in \cite{G} under number (3.49). \end{proof}

\medskip

We go on to Theorem \ref{ext}.

\medskip

\begin{proof*}{Proof of Theorem \ref{ext}} By Theorem \ref{aux} and Lemma \ref{exl}, \begin{align*} \mathcal{A}(a, b) &= \left[\binom{x + n - 1}{n} \bigcirc \mathcal{F}_b(x)\right] + \left[\binom{x + n}{n} \bigcirc \mathcal{F}_b(x - 1)\right]\\ &= \sum_i (-1)^i\binom{1}{i}\mathcal{F}_b(n - i) + \sum_i (-1)^i\binom{0}{i}\mathcal{F}_b(n - i - 1)\\ &= [\mathcal{F}_b(n) - \mathcal{F}_b(n - 1)] + \mathcal{F}_b(n - 1)\\ &= \mathcal{F}_b(n)\\ &= \prod_{i = 1}^s \binom{n + b_i}{b_i}, \end{align*} as needed. \end{proof*}

\medskip

To be fair, our derivation from Theorem \ref{aux} is hardly the simplest way to establish Theorem \ref{ext}. For completeness, let us also sketch one short combinatorial argument.

\medskip

\begin{proof*}{Second proof of Theorem \ref{ext}} For each $i$ with $1 \le i \le s$, let $\beta_i = q_i^{b_i}$ and let $\mathcal{S}_i$ be any canonical series of multiplications transforming $K_{1, 1}$ into $K_{a, \beta_i}$. It is straightforward to see that there are exactly $\binom{n + b_i}{b_i}$ such series.

Out of $\mathcal{S}_1$, $\mathcal{S}_2$, \ldots, $\mathcal{S}_s$, we construct one single canonical series $\mathcal{S}$ of multiplications transforming $K_{1, 1}$ into $K_{a, b}$, as follows.

For each $i$ with $1 \le i \le s$ and each $j$ with $0 \le j \le n$, check if $\mathcal{S}_i$ contains a multiplication on the right which transforms $K_{p^j, u'}$ into $K_{p^j, u''}$ for some $u'$ and $u''$. (Then $u'$ and $u''$ will necessarily be distinct powers of $q_i$.) When it does, set $u_{i, j} = u''/u'$. Otherwise, when it does not, set $u_{i, j} = 1$. Then also let $U_j = u_{1, j}u_{2, j} \cdots u_{s, j}$.

For all $j$ with $0 \le j \le n$ such that $U_j \ge 2$, let $\mathcal{S}$ contain a multiplication on the right by $U_j$ at a moment when the left part of the graph is of size $p^j$. Furthermore, let this account for all multiplications on the right in $\mathcal{S}$. Then the multiplications on the left in $\mathcal{S}$ are determined uniquely as well.

It is straightforward to verify that the mapping $\mathcal{S}_1$, $\mathcal{S}_2$, \ldots, $\mathcal{S}_s \to \mathcal{S}$ we just described is in fact a bijection between all sequences of the form $\mathcal{S}_1$, $\mathcal{S}_2$, \ldots, $\mathcal{S}_s$ and all canonical series of multiplications transforming $K_{1, 1}$ into $K_{a, b}$. The claim follows. \end{proof*}

\medskip

One corollary of Theorem \ref{ext} is that $\mathcal{A}(\{n\}, \{n\}) = \binom{2n}{n}$. Then, by the well-known binomial coefficient identity \[S(2, n) = \binom{2n}{n}\] (number (3.81) in \cite{G}), we conclude that in fact also $\mathcal{A}(\{n\}, \{n\}) = S(2, n)$.

This looks strikingly similar to the identity $\mathcal{A}(\{n, n\}, \{n, n\}) = S(4, n)$ of Theorem~\ref{main}. However, the resemblance appears to be purely coincidental. The author is not aware of any such formula for $\mathcal{A}(\{n, n, n\}, \{n, n, n\})$.

\section{The Proof} \label{proof}

For our proof of Theorem \ref{main}, we follow the same overall strategy as in Section \ref{exp}. This time around, however, the details will be significantly more difficult to fill in.

By Theorem \ref{aux}, \[\mathcal{A}(\{n, n\}, \{n, n\}) = 2\left[\binom{x + n - 1}{n}^2 \bigcirc \binom{x + n}{n}^2\right].\]

The key insight required to evaluate the right-hand side is that we can express $\binom{x + n}{n}^2 \bigcirc P(x)$ neatly in terms of the values of $P$. The exact expression is given by the following lemma.

\medskip

\begin{lemma} \label{square} For all polynomials $P$ with real coefficients, \[\binom{x + n}{n}^2 \bigcirc P(x) = \sum_i (-1)^i\binom{n}{i - n}\binom{i}{n}P(i).\] \end{lemma}

\medskip

We are about to encounter some more complicated binomial coefficient identities which cannot be found in \cite{G}. For them, we need the method of \emph{creative telescoping}.

This method was developed by Doron Zeilberger building upon earlier algorithms of Sister Mary Celine Fasenmyer and Bill Gosper. For a detailed discussion, we point readers to \cite{PWZ}. Here, we give only a quick summary of some key points, and only in the generality required for our purposes.

Let $f(i, n)$ be some concrete expression of $i$ and $n$ such that, for each fixed nonnegative integer $n$, $f(i, n)$ vanishes for all sufficiently large nonnegative integers $i$. Then we can define \[F(n) = \sum_i f(i, n).\]

Suppose that we have somehow managed to find another concrete expression $g(i, n)$ of $i$ and $n$ such that: (a) $g(0, n) = 0$ for all nonnegative integers $n$; (b) For each fixed nonnegative integer $n$, $g(i, n)$ vanishes for all sufficiently large nonnegative integers $i$; and (c) There are concrete polynomials $c_0(n)$, $c_1(n)$, \ldots, $c_k(n)$, all depending only on $n$ and not on $i$, with \[\sum_{j = 0}^k c_j(n)f(i, n + j) = g(i + 1, n) - g(i, n).\]

Then we can conclude immediately that \begin{align*} c_0(n)F(n) + c_1(n)F(n + 1) + \cdots + c_k(n)F(n + k) &= \sum_{j = 0}^k \left[c_j(n)\sum_i f(i, n + j)\right]\\ &= \sum_i \sum_{j = 0}^k c_j(n)f(i, n + j)\\ &= \sum_i [g(i + 1, n) - g(i, n)]\\ &= 0. \end{align*}

Or, in other words, $F$ satisfies a recurrence of order $k$ with coefficients $c_0(n)$, $c_1(n)$, \ldots, $c_k(n)$. The expression $g$ is then called a \emph{certificate} for that recurrence.

The method of creative telescoping allows us to find such a certificate $g$ whenever $f$ is nice in some precise technical sense. The full definition is in \cite{PWZ}, and we do not reproduce it here. What matters is that all four expressions $\varphi$, $\psi$, $\tau$, and $\sigma$ in the proofs of Lemma \ref{square} and Theorem \ref{main} are indeed nice in that way.

The task of discovering the certificate $g$ is usually arduous and best left to machines. Once it is found, however, the verification that it and its associated coefficients $c_0(n)$, $c_1(n)$, \ldots, $c_k(n)$ do indeed work as stated becomes entirely routine. Thus our proofs in this section are still human-friendly, at least in principle. Checking all four of our certificates by hand would probably not be too enjoyable, but certainly it can be done.

There is one important caveat, though. There might be some $n$ such that $g(i, n)$ is not well-defined for all nonnegative integers $i$. Then for those particular values of $n$ we cannot claim that $F$ satisfies our recurrence.

We return to Lemma \ref{square}.

\medskip

\begin{proof*}{Proof of Lemma \ref{square}} Just as with Lemma \ref{exl}, it suffices to consider only the special case when $P(x)$ is of the form $\binom{x}{k}$. Then \[\binom{x}{k} \bigcirc \binom{x + n}{n}^2 = \mathcal{C}\left(k, \binom{x + n}{n}^2\right).\]

On the other hand, by formula ($\ast$) furthermore \[\mathcal{C}\left(k, \binom{x + n}{n}^2\right) = \sum_i (-1)^{k + i}\binom{k}{i}\binom{n + i}{n}^2.\]

With this, we are only left to verify the binomial coefficient identity \[\sum_i (-1)^{k + i}\binom{k}{i}\binom{n + i}{n}^2 = \sum_i (-1)^i\binom{n}{i - n}\binom{i}{n}\binom{i}{k}. \tag{\textbf{A}}\]

Treating $k$ as a fixed parameter, let \[\varphi(i, n) = (-1)^{k + i}\binom{k}{i}\binom{n + i}{n}^2\] and \[\Phi(n) = \sum_i \varphi(i, n).\]

Similarly, let \[\psi(i, n) = (-1)^i\binom{n}{i - n}\binom{i}{n}\binom{i}{k}\] and \[\Psi(n) = \sum_i \psi(i, n).\]

Thus identity (\textbf{A}) takes on the form $\Phi(n) = \Psi(n)$.

By the method of creative telescoping, we find that $\Phi$ satisfies the second-order recurrence \begin{align*} (2n - k + 3)(2n - k + 4)\Phi(n + 2) &= -(n + 1)^2\Phi(n) + {}\\ &\phantom{{} = {}} [5n^2 + (2k + 16)n + 3k + 13]\Phi(n + 1) \end{align*} with certificate \[\frac{\varphi_\text{I}(i, n)}{\varphi_\text{II}(n)} \cdot \varphi(i, n),\] where \[\begin{gathered} \varphi_\text{I}(i, n) = i^3[6n^3 + 7in^2 + 2i^2n - (k - 27)n^2 - {}\\ (2k - 22)in - (k - 3)i^2 - (2k - 40)n - (2k - 18)i - k + 19] \end{gathered}\] and \[\varphi_\text{II}(n) = (n + 1)^2(n + 2)^2.\]

Then we find that $\Psi$ satisfies an identical recurrence as well, with certificate \[\frac{\psi_\text{I}(i, n)}{\psi_\text{II}(n)} \cdot \psi_\text{III}(i, n)\] where \[\begin{gathered} \psi_\text{I}(i, n) = -(i - k)(n - i)^2[12n^3 - 17in^2 + 6i^2n - (k - 55)n^2 + {}\\ (2k - 54)in - (k - 10)i^2 - (2k - 82)n + (2k - 43)i - k + 39],\\ \psi_\text{II}(n) = (n + 1)(n + 2)(n + 3)(n + 4), \end{gathered}\] and \[\psi_\text{III}(n) = (-1)^i\binom{n + 4}{i - n}\binom{i}{n}\binom{i}{k}.\]

Since $\varphi_\text{II}(n)$ and $\psi_\text{II}(n)$ are nonzero when $n \ge 0$, both recurrences hold for all nonnegative integers $n$.

We proceed to check that $\Phi$ and $\Psi$ satisfy the same initial conditions, too. Before we can complete the proof from there, however, there is one more subtlety we must take into account.

Observe that the leading coefficient $(2n - k + 3)(2n - k + 4)$ in the two recurrences vanishes when $n = \lceil k/2 \rceil - 2$. Thus, as far as our recurrences are concerned, anything at all could happen with $\Phi(\lceil k/2 \rceil)$ and $\Psi(\lceil k/2 \rceil)$. Consequently, we must examine these values of $\Phi$ and $\Psi$ manually as well.

For the calculations, we revert back to the original form $\mathcal{C}\left(k, \binom{x + n}{n}^2\right)$ of $\Phi(n)$. By contrast, $\Psi(n)$ does not require special treatment because the number of nonzero terms in it is small in the cases we must consider.

By direct computation, with $n = 0$ we get $\Phi(0) = \Psi(0) = 1$ when $k = 0$ and $\Phi(0) = \Psi(0) = 0$ for all $k \ge 1$. Then, with $n = 1$, we obtain $\Phi(1) = \Psi(1) = 1$ when $k = 0$, $\Phi(1) = \Psi(1) = 3$ when $k = 1$, $\Phi(1) = \Psi(1) = 2$ when $k = 2$, and $\Phi(1) = \Psi(1) = 0$ for all $k \ge 3$. This settles the initial conditions.

With $n = \lceil k/2 \rceil$, we consider two cases based on the parity of $k$.

When $k = 2\ell$ is even, we compare the coefficients before $x^{2\ell}$ on both sides of \[\binom{x + \ell}{\ell}^2 = \sum_i \mathcal{C}\left(i, \binom{x + \ell}{\ell}^2\right)\binom{x}{i}\] to see that $\Phi(\ell) = \mathcal{C}\left(2\ell, \binom{x + \ell}{\ell}^2\right) = \binom{2\ell}{\ell}$. This is also the unique nonzero summand of $\Psi(\ell)$.

Finally, when $k = 2\ell - 1$ is odd, in the same identity we furthermore compare the coefficients on both sides before $x^{2\ell - 1}$. We already know the value of $\mathcal{C}\left(2\ell, \binom{x + \ell}{\ell}^2\right)$ from the previous case, and so in the current case we find that $\Phi(\ell) = \mathcal{C}\left(2\ell - 1, \binom{x + \ell}{\ell}^2\right) = \frac{3}{2}\ell\binom{2\ell}{\ell}$. On the other hand, the only nonzero summands of $\Psi(\ell)$ are now $-\frac{1}{2}\ell\binom{2\ell}{\ell}$ and $2\ell\binom{2\ell}{\ell}$, and we are done. \end{proof*}

\medskip

We are ready to establish Theorem \ref{main}.

\medskip

\begin{proof*}{Proof of Theorem \ref{main}} By Theorem \ref{aux} and Lemma \ref{square}, \begin{align*} \mathcal{A}(\{n, n\}, \{n, n\}) &= 2\left[\binom{x + n - 1}{n}^2 \bigcirc \binom{x + n}{n}^2\right]\\ &= 2\sum_i (-1)^i\binom{n}{i - n}\binom{i}{n}\binom{n + i - 1}{n}^2. \end{align*}

With this, all that is left is to verify the binomial coefficient identity \[2\sum_i (-1)^i\binom{n}{i - n}\binom{i}{n}\binom{n + i - 1}{n}^2 = \sum_i (-1)^{n + i}\binom{2n}{i}^4. \tag{\textbf{B}}\]

Notice that we only need to verify (\textbf{B}) with $n \ge 1$ because to us $n$ is an exponent in the prime factorisations of two positive integers. This is rather fortunate, since (\textbf{B}) is actually false with $n = 0$. On the other hand, when we formally substitute $n = 0$ in the statement of Theorem \ref{main}, we turn out to obtain a correct numerical identity anyway. These two facts might seem at first to contradict each other, but in reality they do not as Theorem \ref{aux} does not hold when $n = 0$ and $a = b = 1$, either.

Back to the proof. Let \[\tau(i, n) = 2 \cdot (-1)^i\binom{n}{i - n}\binom{i}{n}\binom{n + i - 1}{n}^2\] and \[T(n) = \sum_i \tau(i, n).\]

Then also let \[\sigma(i, n) = (-1)^{n + i}\binom{2n}{i}^4,\] and of course we already have the notation $S(4, n)$ for $\sum_i \sigma(i, n)$.

Thus identity (\textbf{B}) takes on the form $T(n) = S(4, n)$.

By the method of creative telescoping, we learn that $T(n)$ satisfies the second-order recurrence \[\begin{gathered} (n + 2)^3(2n + 3)(48n^2 + 66n + 23)T(n + 2) =\\ -4(n + 1)(2n + 1)^3(48n^2 + 162n + 137)T(n) + {}\\ (13056n^6 + 96288n^5 + 289600n^4 + 453428n^3 + 388698n^2 + 172598n + 31030)T(n + 1) \end{gathered}\] with certificate \[\frac{\tau_\text{I}(i, n)}{\tau_\text{II}(n)} \cdot \tau_\text{III}(i, n),\] where \[\begin{gathered} \parbox{0.975\textwidth}{\centering $\scriptscriptstyle \tau_\text{I}(i, n) = (i - 1)(n - i)^2(27552in^{11} - 58176i^2n^{10} + 16032i^3n^9 + 26496i^4n^8 - 12960i^5n^7 + 192i^6n^6 + 96i^7n^5 - 27552n^{11} + 367020in^{10} - 593880i^2n^9 + 69036i^3n^8 + 265008i^4n^7 - 94668i^5n^6 + 936i^6n^5 + 372i^7n^4 - 322620n^{10} + 2113028in^9 - 2572824i^2n^8 - 117016i^3n^7 + 1097868i^4n^6 - 281844i^5n^5 + 1832i^6n^4 + 520i^7n^3 - 1668938n^9 + 6931571in^8 - 6117036i^2n^7 - 1341685i^3n^6 + 2444730i^4n^5 - 436527i^5n^4 + 1780i^6n^3 + 313i^7n^2 - 5021232n^8 + 14336104in^7 - 8555131i^2n^6 - 3650987i^3n^5 + 3159910i^4n^4 - 369970i^5n^3 + 853i^6n^2 + 69i^7n - 9730810n^7 + 19484963in^6 - 6849587i^2n^5 - 5027866i^3n^4 + 2369168i^4n^3 - 162435i^5n^2 + 161i^6n - 12703524n^6 + 17553984in^5 - 2522009i^2n^4 - 3836537i^3n^3 + 953687i^4n^2 - 28891i^5n - 11345120n^5 + 10298636in^4 + 238675i^2n^3 - 1542942i^3n^2 + 159195i^4n - 6891324n^4 + 3757590in^3 + 498924i^2n^2 - 254850i^3n - 2772936n^3 + 789568in^2 + 110216i^2n + 274i^3 - 701016n^2 + 90624in - 2466i^2 - 100936n + 7124i - 6576)$,}\\ \tau_\text{II}(n) = n^2(n + 1)^4(n + 2)(n + 3)(n + 4), \end{gathered}\] and \[\tau_\text{III}(i, n) = (-1)^i\binom{n + 4}{i - n}\binom{i - 1}{n - 1}\binom{n + i - 1}{n}^2.\]

Similarly, we learn that $S(4, n)$ satisfies an identical recurrence as well, with certificate \[\frac{\sigma_\text{I}(i, n)}{\sigma_\text{II}(n)} \cdot \sigma_\text{III}(i, n)\] where \[\begin{gathered} \parbox{0.975\textwidth}{\centering $\scriptscriptstyle \sigma_\text{I}(i, n) = -i^4(n + 1)(7520256n^{14} - 17989632in^{13} + 23126016i^2n^{12} - 24846336i^3n^{11} + 24385536i^4n^{10} - 19193856i^5n^9 + 11046912i^6n^8 - 4521984i^7n^7 + 1303296i^8n^6 - 259584i^9n^5 + 34176i^{10}n^4 - 2688i^{11}n^3 + 96i^{12}n^2 + 143788032n^{13} - 318025728in^{12} + 388985856i^2n^{11} - 405015552i^3n^{10} + 373158144i^4n^9 - 265193472i^5n^8 + 134661120i^6n^7 - 47858688i^7n^6 + 11779104i^8n^5 - 1960896i^9n^4 + 209232i^{10}n^3 - 12720i^{11}n^2 + 324i^{12}n + 1265058816n^{12} - 2581954560in^{11} + 3006342656i^2n^{10} - 3007621632i^3n^9 + 2558678400i^4n^8 - 1615724416i^5n^7 + 711926848i^6n^6 - 215268736i^7n^5 + 44025040i^8n^4 - 5886640i^9n^3 + 477776i^{10}n^2 - 19984i^{11}n + 274i^{12} + 6788332032n^{11} - 12758850560in^{10} + 14127239936i^2n^9 - 13411286784i^3n^8 + 10337382784i^4n^7 - 5692218752i^5n^6 + 2130662560i^6n^5 - 533110848i^7n^4 + 87037664i^8n^3 - 8771336i^9n^2 + 481812i^{10}n - 10412i^{11} + 24823278336n^{10} - 42863544832in^9 + 44967645824i^2n^8 - 39831836800i^3n^7 + 27212428608i^4n^6 - 12767425216i^5n^5 + 3945292272i^6n^4 - 784469984i^7n^3 + 95918420i^8n^2 - 6480836i^9n + 180840i^{10} + 65453976704n^9 - 103508180224in^8 + 102114889152i^2n^7 - 82578817216i^3n^6 + 48702032112i^4n^5 - 18889181056i^5n^4 + 4624614776i^6n^3 - 685321360i^7n^2 + 55815034i^8n - 1897450i^9 + 128392568000n^8 - 185069632256in^7 + 169489802464i^2n^6 - 121699208800i^3n^5 + 59929351944i^4n^4 - 18415238344i^5n^3 + 3348204704i^6n^2 - 328803152i^7n + 13383804i^8 + 190444624224n^7 - 248328030784in^6 + 206880550320i^2n^5 - 127231907728i^3n^4 + 49997564976i^4n^3 - 11395733056i^5n^2 + 1367528712i^6n - 66762840i^7 + 214818976400n^6 - 250240701600in^5 + 183946773576i^2n^4 - 92284271928i^3n^3 + 27026916916i^4n^2 - 4057141508i^5n + 240983000i^6 + 183541567664n^5 - 187075403264in^4 + 115919125400i^2n^3 - 44138004584i^3n^2 + 8536032158i^4n - 632393370i^5 + 117043258224n^4 - 100797821152in^3 + 49016821976i^2n^2 - 12503420152i^3n + 1194417234i^4 + 54070294352n^3 - 37021253024in^2 + 12453372392i^2n - 1586183720i^3 + 17118733664n^2 - 8290920960in + 1433654960i^2 + 3325900992n - 852822144i + 299090304)$,}\\ \sigma_\text{II}(n) = (2n + 1)(2n + 2)^4(2n + 3)^4(2n + 4)^4, \end{gathered}\] and \[\sigma_\text{III}(i, n) = (-1)^{n + i}\binom{2n + 4}{i}^4.\]

Since $\sigma_\text{II}(n)$ is nonzero when $n \ge 0$, the latter recurrence holds for all nonnegative integers $n$. By contrast, $\tau_\text{II}(0) = 0$, and indeed the former recurrence breaks down with $n = 0$. Luckily for us, though, $\tau_\text{II}(n)$ is nonzero when $n \ge 1$. Therefore, our recurrence for $T$ does hold for all positive integers $n$.

Observe also that the leading coefficients of the two recurrences are nonzero for all nonnegative integers $n$, and so this time around we do not encounter any complications similar to the ones in the proof of Lemma \ref{square}.

That only leaves the initial conditions to work through. By direct computation, we see that $T(1) = S(4, 1) = 14$ and $T(2) = S(4, 2) = 786$, as needed. Our proof of Theorem~\ref{main} is complete. \end{proof*}

\section{Acknowledgements} \label{ack}

I would like to thank Professor Donald Knuth for introducing me to the problem of calculating $\mathcal{A}(a, b)$. I would also like to thank Professor Manuel Kauers for helping me apply the method of creative telescoping to the binomial coefficient identities in the proofs of Lemma \ref{square} and Theorem \ref{main}.

\end{document}